\documentclass[a4paper, 12pt]{article}
\usepackage{amssymb, amsmath, amsfonts}
\usepackage{amsthm,amscd,amsfonts,latexsym,color,epsfig}
\begin{document}

\begin{center}
{\bf \Large On a non-local problem for mixed parabolic-hyperbolic type equation with non-smooth line of type changing}\\

\medskip

\textbf{Karimov E.T., Rakhmatullaeva N.A.}\\
E-mail: erkinjon@gmail.com, rakhmatullaeva@mail.ru\\
\smallskip
\emph{Institute of Mathematics, National University of Uzbekistan named after Mirzo Ulughbek (Tashkent, Uzbekistan),}\\
\emph{Tashkent State Technical University named after A.R.Beruni (Tashkent, Uzbekistan)}\\

\end{center}

\bigskip


\bigskip

\textbf{Abstract. }In the present article we investigate a boundary problem with non-local conditions for mixed parabolic-hyperbolic type equation with three lines of type changing. Considered mixed domain contains a rectangle as a parabolic part and three domains bounded by smooth curves and by type-changing lines as a hyperbolic part of the mixed domain. We prove the uniqueness applying energy integral method. The proof of the existence will be done by reducing the original problem into the system of the second kind Volterra integral equations.

\bigskip

\section{Introduction}

Theory of mixed type equations is one of the main parts of the general theory of partial differential equations (PDEs). First fundamental works on this theory were done by Francesco Tricomi [1], S.Gellerstedt [2], F.I.Frankl [3], C.Z.Morawetz [4], M.A.Lavrent'ev and A.S.Bitsadze [5] and etc.

Due to many applications in gas and aero dynamics, mechanics this direction was very rapidly developed. Nowadays, this theory has many branches due to the usage of various methods of mathematical and functional analysis, topological methods and the method of Fractional Calculus.

Omitting huge amount of works, related to the studying mixed type equations, we only note some recent works on local and non-local boundary problems for parabolic-hyperbolic equations [6-10].

Regarding the investigations mixed parabolic-hyperbolic equations with non-smooth lines of type changing we note works [11-14].

In the present paper we study non-local problem for parabolic-hyperbolic equation with three lines of type changing in a special mixed domain, hyperbolic parts of which bounded by smooth curves and by type-changing lines. Under the certain assumptions on these curves, we find conditions to parameters, participated in non-local conditions and we require definite regularity from the given function $f(x,y)$.

\section{Formulation of a problem}

Consider an equation
$$
Lu=f(x,y), \eqno (1)
$$
in a domain $\Omega =\Omega_0\cup\Omega_i\cup AB\cup BC\cup AD$,  where
$$
Lu=\left\{ \begin{array}{l}
u_{xx}-u_y,\hfill (x,y)\in \Omega_0, \\
u_{xx}-u_{yy},\hfill (x,y)\in \Omega_i\,(i=\overline{1,3}) \\
\end{array}
\right.
$$
$\Omega_0$ is a domain bounded by segments $AB, BC, CD, DA$, of straight lines $y=0, x=1, y=1$ and $x=0$ respectively;\\
$\Omega_1$ is a domain bounded by segment $AB$ and smooth curve $\gamma_1:\,y=-\gamma_1(x)$, $\gamma_1(0)=\gamma_1(1)=0$, located inside of the characteristic triangle $0\leq x+y\leq x-y\leq 1$;\\
$\Omega_2$ is a domain bounded by segment $AD$ and smooth curve $\gamma_2:\,x=-\gamma_2(y)$, $\gamma_2(0)=\gamma_2(1)=0$, located inside of the characteristic triangle $0\leq y+x\leq y-x\leq 1$;\\
$\Omega_3$ is a domain bounded by segment $BC$ and smooth curve $\gamma_3:\,x=-\gamma_3(y)$, $\gamma_3(0)=\gamma_3(1)=1$, located inside of the characteristic triangle $1\leq x-y\leq x+y\leq 2$.

Regarding the curves $\gamma_i(t)\,(i=\overline{1,3})$  we suppose that they are twice continuously differentiable and $t\pm\gamma_i(t)\,(0\leq t\leq 1, i=\overline{1,3})$  are monotonically increase.

$\theta_1(t),\theta_2(t),\theta_3(t)\left[\theta_1^*(t),\theta_2^*(t),\theta_3^*(t) \right]$ are affixes of points of intersection of the curves $\gamma_i(t)\,(i=\overline{1,3})$ and characteristics $x-y=t, y-x=t, x+y=1+t$ $\left[x+y=t, x-y=t, y-x=1+t\right]$ of the Eq.(1), respectively.

We formulate the following problem for Eq. (1).

\textbf{Problem.} To find a regular solution of Eq. (1), satisfying nonlocal conditions
$$
\left[u_x-u_y\right]\left(\theta_1(t)\right)=\sigma_1\left[u_x+u_y\right]\left(\theta_1^*(t) \right),0\le t\le 1, \eqno (2)
$$
$$
\left[u_x-u_y\right]\left(\theta_2(t)\right)=\sigma_2\left[u_x+u_y\right]\left(\theta_2^*(t) \right),0\le t\le 1, \eqno (3)
$$
$$
\left[u_x+u_y\right]\left(\theta_3(t)\right)=\sigma_3\left[u_x-u_y\right]\left(\theta_3^*(t) \right),0\le t\le 1, \eqno (4)
$$
$$
u(A)=u(B)=0. \eqno (5)
$$
Here $\sigma_1, \sigma_2, \sigma_3$ are arbitrary real numbers.

We call as a regular solution of the problem in the domain $\Omega$ function $u(x,y)\in W$, where $W=\left\{ u:u(x,y)\in C\left( \overline{\Omega } \right)\cap C^1\left( \bar{\Omega } \right)\cap C_{x,y}^{2,1}\left( \Omega_0 \right)\cap C^2\left( \Omega_i \right),i=\overline{1,3} \right\}$
satisfying Eq.(1) in domains $\Omega_i\, (i=\overline{0,3})$.

\section{Main result}

\textbf{Theorem.} If conditions $\sigma_2,\sigma_3\in [-1,1]$, $f(x,y)\in C^2(\Omega)$  are fulfilled, then the problem has unique regular solution.

\verb"Proof:" First, we deduce main functional relations.

Solution of the problem in $\Omega_i\, (i=\overline {1,3})$ can be represented by the D'Alembert's formula [15]:
$$
u(\xi ,\eta )=\frac{1}{2}\left[ \tau_1^-(\xi )+\tau_1^-(\eta )-\int\limits_{\xi }^{\eta }{\nu_1^-(t)dt} \right]-\int\limits_{\xi }^{\eta }d\xi_1\int\limits_{\xi_1}^{\eta}{f_1\left(\xi_1,\eta_1\right)d\eta_1}, \eqno (6)
$$
$$
u(\xi ,\eta )=\frac{1}{2}\left[ \tau_2^-(\xi )+\tau_2^-(-\eta )-\int\limits_{\xi }^{-\eta }{\nu_2^-(t)dt} \right]-\int\limits_{\xi }^{-\eta }d\xi_1\int\limits_{\xi_1}^{-\eta}{f_1\left(\xi_1,\eta_1\right)d\eta_1}, \eqno (7)
$$
$$
u(\xi ,\eta )=\frac{1}{2}\left[ \tau_3^+(\xi-1 )+\tau_3^+(1-\eta )-\int\limits_{\xi -1}^{1-\eta }{\nu_3^+(t)dt} \right]-\int\limits_{\xi-1 }^{1-\eta }d\xi_1\int\limits_{\xi_1}^{1-\eta}{f_1\left(\xi_1,\eta_1\right)d\eta_1}, \eqno (8)
$$
where $\xi=x+y, \eta=x-y$, $4f_1(\xi,\eta)=f\left(\frac{\xi+\eta}{2}, \frac{\xi-\eta}{2}\right)$.

By virtue of conditions to $\gamma_i\,(i=\overline{1,3})$, equation of the curve represented as $\xi =\rho (\eta )$  and $\eta =\upsilon (\xi )$ such that $\rho \left( \upsilon (\xi ) \right)=\xi$.

Since
$$
\theta_1\left( \frac{x-\gamma_1(x)+t}{2};\frac{x-\gamma_1(x)-t}{2} \right);\,\,
\theta_1^*\left(\frac{x+\gamma_1(x)+t}{2};-\frac{x+\gamma_1(x)+t}{2} \right);
$$
$$
\theta_2\left( \frac{y-\gamma_2(y)-t}{2};\frac{y-\gamma_2(y)+t}{2} \right);\,\,
\theta_2^*\left( -\frac{y+\gamma_2(y)-t}{2};\frac{y+\gamma_2(y)+t}{2} \right);
$$
$$
\theta_3\left(- \frac{y+\gamma_3(y)-1-t}{2};\frac{y+\gamma_3(y)+1+t}{2} \right);\,\,
\theta_3^*\left( \frac{y-\gamma_3(y)+1-t}{2};\frac{y-\gamma_3(y)-1+t}{2} \right),
$$
from (2) and (6), (3) and (7), (4) and (8) we deduce main functional relations on the line of type changing of Eq.(1), respectively:
$$
(1-\sigma_1){\tau_1^-}'(x)-(1+\sigma_1)\nu_1^-(x)=A_1(x),\,\,\,0<x<1, \eqno (9)
$$
$$
(1+\sigma_2){\tau_2^-}'(y)+(1-\sigma_2)\nu_2^-(y)=A_2(y),\,\,\,0<y<1, \eqno (10)
$$
$$
(1+\sigma_3){\tau_3^+}'(y)+(1-\sigma_3)\nu_3^+(y)=A_3(y),\,\,\,0<y<1, \eqno (11)
$$
where
$$
A_1(x)=2\sigma_1\int\limits_x^{\upsilon (x)}{f_1\left(x;\eta_1\right)d\eta_1}+2\int\limits_{\rho (x)}^x{f_1\left(\xi_1;x\right)d\xi_1},
$$
$$
A_2(y)=-2\sigma_2\int\limits_y^{\upsilon (y)}{f_1\left(x;\eta_1\right)d\eta_1}+2\int\limits_{\rho (y)}^y{f_1\left(\xi_1;x\right)d\xi_1},
$$
$$
A_3(y)=2\sigma_3\int\limits_y^{\upsilon (y)}{f_1\left(x;\eta_1\right)d\eta_1}-2\int\limits_{\rho (y)}^y{f_1\left(\xi_1;x\right)d\xi_1}.
$$

\subsection{The uniqueness of the solution of the problem}

Consider the case, when $\left|\sigma_i \right|\ne 1$ .

In the domain $\Omega_0$  we have the equality
$$
\begin{array}{l}
\displaystyle{\iint\limits_{\Omega_0}{u_x^2(x,y)dxdy}+\int\limits_0^1{\tau_2^+(y)\nu_2^+(y)dy}-\int\limits_0^1{\tau_3^-(y)\nu_3^-(y)dy}+}\hfill\\
\displaystyle{+\frac{1}{2}\int\limits_0^1{u^2(x,1)dx}-\int\limits_0^1{\left[ \tau_1^+(x) \right]^2dx=0.}}\hfill\\
\end{array}
\eqno (12)
$$

In order to prove the uniqueness, first we prove that $u(x,\pm 0)={\tau_1}^{\pm }(x)=0$. Passing to the limit in $\Omega_0$, at $y\to +0$, from the equation $u_{xx}-u_y=0$, we get ${\tau_1^+}''(x)=\nu_1^+(x)$ and substituting it into the integral $I_1=\int\limits_0^1{\tau_1^+(x)\nu_1^+(x)dx}$, taking condition (5) into account we have
$$
I_1=\int\limits_0^1{\tau_1^+(x){\tau_1^+}''(x)dx}=-\int\limits_0^1{\left(\left(\tau_1^+(x)\right)'\right)^2dx}.\eqno (13)
$$
Now considering (9) we obtain
$$
I_1=\frac{1-\sigma_1}{1+\sigma_1}\int\limits_0^1{\tau_1^+(x){\tau_1^+}'(x)dx}=\left. \frac{1-\sigma_1}{2(1+\sigma_1)}\left(\tau_1^+(x) \right)^2 \right|_0^1=0. \eqno (14)
$$
From (13) and (14) it follows that $\tau_1^{\pm }(x)\equiv 0$.

Now we prove that
$$
I_2=\int\limits_0^1{\tau_2^+(y)\nu_2^+(y)dy\ge 0},\,\,I_3=\int\limits_0^1{\tau_3^+(y)\nu_3^+(y)dy\le 0}.
$$
Taking (10) and (11) into account we get respectively the followings:
$$
I_2=\int\limits_0^1{\tau_2^+(y)\nu_2^+(y)dy}=\frac{1+\sigma_2}{1-\sigma_2}\int\limits_0^1{\tau_2^+(y){\tau_2^+}'(y)dy}=\frac{1+\sigma_2}{2(1-\sigma_2)}
{\tau_2^+}^2(1),
$$
$$
I_3=\int\limits_0^1{\tau_3^+(y)\nu_3^+(y)dy}=\frac{1+\sigma_3}{1-\sigma_3}\int\limits_0^1{\tau_3^-(y){\tau_3^-}'(y)dy}=-\frac{1+\sigma_3}{2(1-\sigma_3)}
{\tau_3^-}^2(1).
$$
If $\frac{1+\sigma_2}{1-\sigma_2}>0$, $\frac{1+\sigma_3}{1-\sigma_3}>0$, then $I_2\geq 0, I_3\leq 0$.

Considering that $I_2\geq 0, I_3\leq 0$ and $\tau_1^+=0$ from the equality (12) we state that $u(x,y)=0$  in $\Omega_0$ . Since $u(x,y)\in C(\overline{\Omega})$, we can conclude, $u(x,y)=0$  in  $\Omega$.

\subsection{The existence of the solution of the problem}

Now we prove the existence of the solution for the problem.

Passing to the limit in $\Omega_0$  at $y\to +0$, from the equation $u_{xx}-u_y=g(x;y)$, taking (9) into account we get
$$
{\tau_1^+}''(x)-\frac{1-\sigma_1}{1+\sigma_1}{\tau_1^+}'(x)=f^*(x),\eqno (15)
$$
where $f^*(x)=f(x,0)-\frac{1}{1+\sigma_1}A_1(x).$
From the condition (5) we have 
$$
\tau_1^+(0)=0,\,\,\tau_1^+(1)=0. \eqno (16)
$$
Solution of Eq. (15) together with conditions (16) can be represented as
$$
\tau_1^+(x)=\frac{1+\sigma_1}{1-\sigma_1}\left[ \int\limits_0^x{\left( e^{\frac{1-\sigma_1}{1+\sigma_1}(x-t)}-1 \right)f^*(t)dt}-\frac{e^{\frac{1-\sigma_1}{1+\sigma_1}x}-1}{e^{\frac{1-\sigma_1}{1+\sigma_1}}-1}\int\limits_0^1{\left( e^{\frac{1-\sigma_1}{1+\sigma_1}(1-t)}-1 \right)f^*(t)dt} \right].
$$

Solution of the first boundary problem for the Eq. (1) in the domain $\Omega_0$ has a form [16]
$$
u(x,y)=\int\limits_0^1{\tau_1^+(x_1)G(x,y;x_1,y_1)dx_1}+\int\limits_0^y{\tau_2^+(y_1)G_{x_1}(x,y;0,y_1)dy_1}-
$$
$$
-\int\limits_0^y{\tau_3^-(y_1)G_{x_1}(x,y;1,y_1)dy_1}-\int\limits_0^1{dx_1
\int\limits_0^y{f(x_1,y_1)G(x,y;x_1,y_1)dy_1.}}
$$
where 
$
G(x,y;x_1,y_1)=\frac{1}{2\sqrt{\pi (y-y_1)}}\,\sum\limits_{n=-\infty }^{\infty }{\left[ e^{-\frac{(x-x_1+2n)^2}{4(y-y_1)}}-e^{-\frac{(x+x_1+2n)^2}{4(y-y_1)}} \right]}
$ is the Green's function of the first boundary problem for the heat equation.

Differentiating (17) once by $x$ and considering (10), (11), introducing notation 
$$
N(x,y;x_1,y_1)=\frac{1}{2\sqrt{\pi (y-y_1)}}\,\sum\limits_{n=-\infty }^{\infty }{\left[ e^{-\frac{(x-x_1+2n)^2}{4(y-y_1)}}+e^{-\frac{(x-x_1+2n)^2}{4(y-y_1)}} \right]},
$$
and taking 
$$
G_{x_{_1}x}(x,y;x_1,y_1)=N_{y_{_1}}(x,y;x_1,y_1),\,\,\,G_x(x,y;x_1,y_1)=-N_{x_{_1}}(x,y;x_1,y_1)
$$ 
into account after the integration by part we obtain the followings:
$$
\frac{1+\sigma_2}{1-\sigma_2}{\tau_2^+}'(y)+\int\limits_0^y{{\tau_2^+}'(y_1)N(0,y;0,y_1)dy_1}-\int\limits_0^y{{\tau_3^-}'(y_1)N(0,y;1,y_1)dy_1}
=E_1(y),
\eqno (18)
$$
$$
\frac{1+\sigma_3}{1-\sigma_3}{\tau_3^-}'(y)-\int\limits_0^y{{\tau_3^-}'(y_1)N(1,y;1,y_1)dy_1}+\int\limits_0^y{{\tau_2^+}'(y_1)N(1,y;0,y_1)dy_1}=
E_2(y),\eqno (19)
$$
where 
$$
E_1(y)=\frac{A_2(y)}{1-\sigma_2}+\int\limits_0^1{{\tau_1^+}'(x_1)N(0,y;x_1,y_1)dx_1}+\int\limits_0^1{dx_1
\int\limits_0^y{f(x_1,y_1)G_x(0,y;x_1,y_1)dy_1}},
$$
$$
E_2(y)=\frac{A_3(y)}{1-\sigma_3}+\int\limits_0^1{{\tau_1^+}'(x_1)N(1,y;x_1,y_1)dx_1}+\int\limits_0^1{dx_1\int\limits_0^y{f(x_1,y_1)G_x(1,y;x_1,y_1)dy_1}}.
$$

Since (18) is the second kind Volterra integral equation regarding the function ${\tau_2^+}'(y)$, we can formally write its solution via resolvent kernel.  Substituting found solution ${\tau_2^+}'(y)$ into (19), we will get second kind Volterra integral equation regarding the function ${\tau_3^-}'(y)$, which is uniquely solvable. After the finding of functions $\tau_i^{\pm }(t)\,\left( i=\overline{1,3} \right)$, we find unknown functions $\nu_i^{\pm }(t)\,\left( i=\overline{1,3} \right)$ by formulas (9)-(11).

Since we have all functions required to write solution of the problem, it can be represented in the domain   by the formula (17), in domains $\Omega_i(i=\overline{1,3})$ by formulas (17) - (19), respectively.

In case, when $\left|\sigma_i\right|=1$ the problem will be divided into 4 problems, because by formulas (9)-(11) we can find ${\tau_i^{\pm }}'(t)$ or $\nu_i^{\pm }(t)$ $\left(i=\overline{1,3} \right)$.

Theorem is proved. 

\smallskip

\centerline{\bf References}
\begin{enumerate}
\item {\it Tricomi F. G.} {Sulle Equazioni Lineari alle derivate Parziali di 2? Ordine, di Tipo Misto. //  Atti Accad. Naz. dei Lincei. 14 (5) 1923, p. 133-247}\\
\item{\it Gellerstedt S.} {Sur un probleme aux limits pour une equation liniare aux derivees du second orde de type mixte. // Dissertation. Uppsala University. 1935} \\
\item {\it Frankl F.I.} {Selected works on gas dynamics. Moscow:Nauka, 1973, 711 p.}\\
\item{\it Morawetz C.Z.} {A weak solution for a system of equations of elliptic-hyperbolic type. // Comm. Pure and Appl. Math. 1958. -V.11.-P. 315-331} \\
\item{\it Lavrent'ev M.A., Bitsadze A.V.} {To the problem of mixed type equations. Dokl. AN SSSR, Vol.70, No 3, 1950, pp.373-376.}\\
\item{\it Berdyshev A.S. and Karimov E.T.} {Some non-local problems for the parabolic-hyperbolic type equation with non-characteristic line of changing type. Cent. Eur. J. Math., Vol. 4, (2006), no. 2, pp. 183-193.}\\
\item{\it Ashyralyev A., Yutsever H.A.} {On difference schemes for hyperbolic-parabolic equations. Proceedings of the int. conf. "Dynamical systems and applications", Antalya, Turkey, 5-10 July, 2004, pp. 136-153.}\\
\item{\it Berdyshev A.S.} {Basis property of the system of radical functions of the non-local boundary-value problem for the parabolic-hyperbolic equation, Dokl. Akad. Nauk 366 (1999), No. 1, pp.7-9 }\\
\item{\it Karimov E.T.} {Some non-local problems for the parabolic-hyperbolic type equation with complex spectral parameter. Mathematische Nachrichten, vol. 281 (2008), no. 7, pp. 959-970.}\\
\item{\it Sabytov K.B.} {To the theory of mixed parabolic-hyperbolic type equations with a spectral parameter, Differentsial'nye Uravneniya 25 (1989), No. 1, pp.117-126}\\
\item{\it Nakhusheva V.A.} {First boundary problem for mixed type equation in a characteristic polygon. Dokl.AMAN, 2012. Vol.14, No 1, pp.58-65.}\\
\item{\it Berdyshev A.S., Rakhmatullaeva N.A.} {Non-local problems for parabolic-hyperbolic equations with deviation from the characteristics and three type-changing lines. Electronic Journal of Differential Equations. Vol. (2011) 2011, No 7, pp.1-6.} \\
\item{\it Eleev V.A., Lesev V.N.} {On two boundary problems for mixed type equations with perpendicular lines of type changing// Vladikavkaz math.journ. 2001. Vol. 3. Vyp.4, pp.9-22.}\\
\item{\it G.D. Tojzhanova and M.A. Sadybekov} {About spectral properties of one analogue of the Tricomi problem for the mixed parabolic-hyperbolic type equation. Izvestija AN KazSSR, ser.phys.-math.nauk, Vol.3 (1989), pp. 48-52.}\\
\item {\it Rassias J.M.} {Lecture Notes on Mixed Type Partial Differential Equations. World Sci., Singapore, 1990.}\\
\item {\it Djuraev T.D., Sopuev A., and Mamajonov M.} {Boundary-value Problems for the Parabolic-hyperbolic Type Equations (Fan, Tashkent, 1986).}\\

\end{enumerate}
\end{document}